\documentclass[reqno,12pt]{amsart}
\usepackage{amssymb}

\newcommand{\N}{\hbox{I\hskip -2pt N}}

\newcommand{\ra}{\rightarrow}

\newtheorem{theo}{Theorem}

\newtheorem{conj}{Conjecture}

\newtheorem{lem}{Lemma}

\def\pro{\noindent {\bf Proof. }}

\begin{document}

\title {Coloring dense graphs via VC-dimension}
\author{Tomasz {\L}uczak}
\address{Faculty of Mathematics and CS
Adam Mickiewicz University
Poznan, Poland {\tt tomasz@amu.edu.pl}}
\thanks{The first author partially supported by the
Foundation for Polish Science.}
\author{St\'ephan Thomass\'e}
\address{Universit\'e Montpellier 2 - CNRS, LIRMM 161 rue Ada, 34392 Montpellier
    Cedex, France
    {\tt thomasse@lirmm.fr}}
\date{January 26, 2010}

\maketitle

\begin{abstract}
The Vapnik-\v{C}ervonenkis dimension is a complexity measure
of set-systems, or {\it hypergraphs}. Its application to graphs
is usually done by considering the sets of neighborhoods
of the vertices (see~\cite{ABK} and~\cite{CEV}), hence providing a set-system. But the graph structure
is lost in the process. The aim of this paper is to introduce the notion of
{\it paired VC-dimension}, a generalization
of VC-dimension to set-systems endowed with a graph structure,
hence a collection of pairs of subsets.

The classical VC-theory is generally
used in combinatorics to bound the transversality
of a hypergraph in terms of its fractional transversality
and its VC-dimension. Similarly, we
bound the chromatic number
in terms of fractional transversality and paired VC-dimension.

This approach turns out to be very useful for
a class of problems raised by Erd{\H o}s and Simonovits~\cite{ES} asking for
$H$-free graphs with minimum degree at least $cn$ and arbitrarily high chromatic number,
where $H$ is a fixed graph and $c$ a positive constant.

We show how the usual VC-dimension
gives a short proof of the fact that triangle-free
graphs with minimum degree at least $n/3$ have bounded
chromatic number, where $n$ is the number of vertices.

Using paired VC-dimension, we prove that if the chromatic number of $H$-free graphs
with minimum degree at least $cn$ is unbounded for
some positive $c$, then it is unbounded for all $c<1/3$.
In other words, one can find $H$-free graphs with unbounded
chromatic number and minimum degree arbitrarily close to $n/3$. These $H$-free graphs are
derived from a construction of Hajnal. The large chromatic number follows
from the Borsuk-Ulam Theorem.
\end{abstract}

\section{Introduction.}

The {\it chromatic number} of a graph $G=(V,E)$, i.e. the minimum number
of parts of a partition of its vertex set into edgeless subsets ({\it stable sets}),
is one of the most studied parameters in graph theory. However, this parameter
cannot be directly interpreted as a measure of the complexity of $G$
since very simple graphs, like {\it cliques} $K_n$ of size $n$ inducing
all possible edges, have chromatic number $n$.

The picture becomes completely different when a graph $G$ has
large chromatic number for non obvious reasons like the
containment of a large clique. For instance, in the case of
triangle-free graphs, achieving high chromatic number is not
a straightforward exercise. Indeed this simple question is the
starting point of several areas like random  and topological
graphs.

Specifically some constructions of triangle-free graphs with high
chromatic number were provided first by Zykov~\cite{Zyk}, and then
by Mycielski~\cite{Myc}. Erd\H{o}s~\cite{ER} proposed a construction based
on random graphs with arbitrarily large girth. Geometric constructions based
on the Borsuk-Ulam theorem provide examples with arbitrarily large odd girth.
However, all graphs constructed in the above way are sparse, i.e. have small
minimum degree with respect to the number of vertices.

In their seminal paper, Erd{\H o}s and Simonovits~\cite{ES} asked for
a bound on the chromatic number of triangle-free graphs with minimum
degree larger than $n/3$. This question was first solved by Thomassen in~\cite{CT1}
where he provided a bound on the chromatic number when the
degree is larger than $(1/3+\varepsilon)n$, for every $\epsilon >0$.
This result was sharpened by {\L}uczak~\cite{TL} who showed that,
up to homomorphism, there are
only finitely
many maximal triangle-free graphs with minimum
degree larger than $(1/3+\varepsilon)n$. The finiteness relies on the partition
provided by the regularity lemma. Finally, Brandt and Thomass\'e~\cite{BT} proved
that all such graphs have chromatic number at most
four, using a complete characterization of the family.

The goal of our paper is to show how to use  Vapnik-\v{C}ervonenkis theory
can be used for the Erd{\H o}s-Simonovits type problems. The three main
result of this paper are the following.

\begin{itemize}
\item We give a new short proof of the existence of a bound
on the chromatic number of triangle-free graphs with minimum
degree larger than $n/3$. This direct consequence of classical
VC-dimension even allows to break the $1/3$ barrier. For instance,
it provides a bound for minimum degree $n/3$ minus a constant.
This was completely out of reach using the previous methods of~ \cite{BT}, \cite{TL}, and~\cite{CT1}.

\item Introducing a new parameter, the {\it paired VC-dimension}, we characterize the
graphs $H$ such that the class of $H$-free graphs (with respect
to homomorphism) has chromatic threshold 0 (i.e. have bounded
chromatic number as soon as the minimum
degree is larger than $cn$, where $c>0$). From this characterization
follows that the chromatic threshold of $H$-free graphs is either
0 or at least 1/3. For instance 0 is achieved by the pentagon,
1/3 by the triangle, but no value in between can be realized.

\item Thomassen~\cite{CT2} recently proved that the class of
pentagon-free graphs (with respect to subgraph) has chromatic
threshold 0. We give a new proof of this result based on
paired VC-dimension. Using our method we construct a wide class of
non-bipatite graphs (which includes graph without a copy
of the Petersen graph)
for which the threshold is also 0.
\end{itemize}

The central theorem of the paper is  Theorem~\ref{main}
which gives an upper bound on the chromatic number of a graph in terms
of its minimum degree and paired VC-dimension. We believe that
Theorem~\ref{main} can be both generalized and sharpened
(the implicit bounds for the chromatic number which follow
from our argument are rather poor). Another key
observation of the paper is the "duality" between VC-dimension,
which provides upper bounds
on the chromatic number of dense graphs, and the Borsuk-Ulam theorem,
which gives constructions achieving lower bounds.

The structure of the paper goes as follows.
We start with a construction of dense
graphs   based on Borsuk graphs,
i.e. graphs living on the
$d$-dimensional sphere where edges link
near-antipodal points. They have chromatic number $d+2$
according to Borsuk-Ulam theorem.
We shall use them later to show
that the chromatic threshold jumps from 0 to 1/3.
Then we state
Erd\H{o}s--Simonovits  problem for general graphs $H$.
The next section contains a short argument which applies
the standard version of VC-dimension to study the
chromatic number of triangle-free graphs. Only then we
introduce the key notion of paired VC-dimension and
prove our main result Theorem~\ref{main} from which
the chromatic thresholds for $H$ (Theorems~\ref{homfree} and
\ref{subfree})
naturally follow. We conclude the paper with some
remarks and open problems.

\section{Borsuk-Hajnal Graphs.}

The first construction of a dense triangle-free graph with high chromatic number
is due to Hajnal. It is based on Kneser graphs, which confer the large chromatic number,
to which is added a bipartite graph of large size, in order to raise
the minimum degree.

Recall that the {\it Kneser graph} $Kn(n,k)$ has vertex-set $2n+k\choose n$, the collection of $n$ elements
subsets of the set $\{1,\dots ,2n+k\}$. Two vertices $x,y$ of $Kn(n,k)$ are joined
by an edge if and only if $x\cap y=\emptyset$.

There exists a canonical coloring of $Kn(n,k)$ into $k+2$ colors: color by 1 the vertices
which contain 1, by 2 the remaining vertices which contain 2, ... until the color $k$.
The vertices not colored up to this point are the $n$-element subsets of $\{k+1,\dots ,2n+k\}$,
hence they induce a matching of $Kn(n,k)$, which is indeed 2-colorable.
In all, $k+2$ colors have been used.
The celebrated theorem of Lov\'asz~\cite{LL} asserts that this coloring is optimal, and thus that $Kn(n,k)$
is $k+2$ chromatic.
This was the starting point of Hajnal's construction (let us remark however,
that not knowing Lov\'asz' result he used more crude lower bound for
$\chi(Kn(n,k)$).

Consider three
integers $k\ll n\ll \ell$, with the property that $2n+k$ divides $\ell$. Form the disjoint
union of $Kn(n,k)$ and a stable set $S$ of size $2\ell$. Divide $S$  into $2n+k$ equal
parts $S_1,\dots ,S_{2n+k}$. Recall that a vertex $x$ of $Kn(n,k)$ corresponds to some $n$-subset of
$2n+k$. Link $x$ to every subset $S_i$ for which $i\in x$. The last step of the construction
consists of adding another stable set $S'$ of size $\ell$ completely joined
to~$S$.

Let us denote this {\it Hajnal graph} by $H(k,n,\ell)$.  It satisfies the following properties:

\begin{itemize}
\item Its chromatic number is at least $k+2$.
\item It does not contain a triangle.
\item Its minimum degree is at least $2\ell n/(2n+k)$, while its total
number of vertices is $3\ell + {2n+k\choose n}$. Hence, since
$k\ll n\ll \ell$, the size of the neighborhood of every vertex is just below a third of the total number
of vertices.
\end{itemize}

The existence of Hajnal graphs directly implies the following results.

\begin{theo}\label{threshtriangle}
For every $\varepsilon>0$ and every integer $k$, there exists a triangle-free graph
on $n$ vertices with minimum degree at least $(1/3-\varepsilon) n$ and chromatic number at least $k$.
\end{theo}

However, for our purpose, we need an alternate construction of Hajnal graphs, which
only differ by the subgraph used to provide the large chromatic number. Here the Kneser graph
is replaced by a Borsuk graph.

Given an integer $d$ and a small $\varepsilon >0$, the {\it Borsuk
graph} $Bor(d,\varepsilon)$ has vertex-set the points of the unit sphere
$S^d$ and edge set the pairs of vertices $xy$ which distance on the sphere
is at least  $\pi-\varepsilon$, where here and below we measure the distance
between two point on the sphere by the angle between them.

This graph has $2^{\aleph _0}$ vertices and Borsuk-Ulam theorem asserts that
its chromatic number is $d+2$ (this fact being easily equivalent  to Borsuk-Ulam theorem).
However the minimum degree of the Borsuk graph is very small. We now increase
it as in Hajnal graph.

Let $\delta >0$ and $x$ be a point of the sphere $S^d$. The {\it $(x,\delta)$-cap} is
the set of vertices of the sphere which distance to $x$ is at most $\pi/2-\delta$.
Now form a graph with parameters $(d,\varepsilon,\delta)$  as follows:

\begin{itemize}
\item start with $B$, a Borsuk graph  $Bor(d,\varepsilon)$;
\item add a stable set $S$ which vertices correspond to all $(x,\delta)$-caps;
\item join $x\in B$ to $c\in S$ if the point $x$ belongs to the cap $c$;
\item add a vertex $v$ joined to $S$.
\end{itemize}

To achieve our goal, we now need to discretise this structure, and add some
weight on the vertices to reach minimum degree close to 1/3.
To this end consider a finite, however dense enough,
set $B'$ of vertices of $B$ and set $S'$ of caps of $S$ still achieving
chromatic number $d+2$ and such that, provided $\delta $ is small enough, every
vertex of $B'$ is joined to nearly half of the vertices of $S'$. We also keep the vertex $v$
to which is given weight $1/3$, while weight $2/3$ is equally distributed on the vertices of $S'$. 
The vertices of $B'$ receive for themselves some negligible weight.

The weighted graph we just constructed is such that the weight of the neighborhood of
every vertex is at least just below $1/3$ and has chromatic number at least $d+2$.
If needed, one can easily turn this weighted graph into a real one by blowing up
every vertex by a stable set of size proportional to its weight. Let us call these
graphs {\it Borsuk-Hajnal} graphs with parameters $d,\varepsilon,\delta$.

Since this graph is also triangle-free, this is an alternate construction for Hajnal
graph. However, they enjoy the following additional property, which we will need
for our main result.

\begin{lem} \label{oddgirth}
In a Borsuk-Hajnal graph with parameter $0<\varepsilon\ll \delta$, where $\delta$ is small enough, every short odd cycle has at least two vertices outside~$B'$.
\end{lem}

\pro The odd girth of a Borsuk graph can be made arbitrarily
large, provided $\varepsilon$ is chosen small enough. Hence no
short odd cycle can be contained in $B'$. This property is also
valid inside the Kneser graph of Hajnal graphs, but there exists
cycles of length five containing four vertices inside the
Kneser graph. Let us prove that this cannot be the case in Borsuk-Hajnal
graphs

To see this, observe that if two vertices $x$ and $y$ of $B'$ are joined in $B'$
by a path of length $2k+1$, then the distance on $S^d$ from $x$ to $y$ is at least
$\pi -(2k+1)\varepsilon$. Hence if $\delta $ is greater than $(2k+1)\varepsilon$,
no $\delta $-cap contain both $x$ and $y$. Consequently,
 every odd cycle longer than $2k+3$
has at least two vertices outside $B'$.~$\square$

~

By the {\it near odd girth} of a Borsuk-Hajnal graph
we mean the minimum length of an odd cycle which has at most one vertex
outside of $B'$. Lemma~\ref{oddgirth} asserts that the near
odd girth can be chosen arbitrarily large.

\section{Chromatic Threshold.}

Let $\mathcal  G$ be a class of (finite) graphs. The
{\it chromatic threshold} of $\mathcal  G$ is the infimum
of the values $c\geq 0$ such that the subclass of
$\mathcal  G$ consisting of the graphs $G$ with minimum
degree at least $c |V(G)|$ have bounded chromatic number.

When $\mathcal  G$ is the class of triangle-free graphs, Theorem~\ref{threshtriangle}
asserts that the chromatic threshold of $\mathcal  G$ is at least
$1/3$.
The study of the chromatic threshold was initiated by Erd\H{o}s and
Simonovits in \cite{ES}, where the problem was posed for
triangle-free graphs and other classes.

The answer for triangle-free graphs was first given by
Thomassen in \cite{CT1} where he proved that the value
was indeed 1/3. This concluded a long list of results.
For instance, Andr\'asfai, Erd\H{o}s and S\'os~\cite{AES} proved
that triangle-free graphs with degree at least $2|V(G)|/5$
are bipartite. Jin~\cite{JIN} showed
that degree at least $10|V(G)|/29$ gives 3-colorable graphs.
These two results being sharp.

Two new proofs of the threshold 1/3 were later given by
{\L}uczak~\cite{TL} using the regularity lemma, and by
Brandt and Thomass\'e~\cite{BT}. This last paper gives
a complete characterization of triangle-free graph with
minimum degree larger than $|V(G)|/3$ and
show that all of them are 4-colorable.

However, the methods used in the above papers could not be extended
to degree at least $|V(G)|/3$ (instead of more than).

The chromatic threshold for $K_r$-free graphs was recently
found by Goddard and Lyle \cite{GoLy}. They used the fact  that,
in principle, the structure of  $K_r$-free graphs inherits
the structure of triangle-free graphs, to show
that the threshold is equal to $(2r-5)/(2r-3)$ in this case.

However, when forbidding other graphs, the problem is wide open.
For instance, a recent
result of Thomassen~\cite{CT2} (answering a question of Erd\H{o}s and
Simonovits) asserts that the chromatic threshold of $C_5$-free graphs
is $0$. In other words, one cannot expect dense graphs without
cycles of length five and with large chromatic number.

Let us introduce two notations. Given a fixed graph $H$ the class
of {\it $H$-hom-free} graphs consists of the graphs without a homomorphic
copy of $H$, while {\it $H$-sub-free} graphs do not contain a subgraph
$H$.

Our main result in this paper is to characterize the graphs $H$
for which the class of $H$-hom-free graphs has chromatic threshold 0.
More specifically, we prove that such an $H$ has this property
if and only if it embeds in all large enough Borsuk-Hajnal Graphs.
In particular, if the chromatic threshold is not zero, since
the minimum degree of Borsuk-Hajnal Graphs is close to 1/3,
the chromatic threshold is at least 1/3.

Moreover our method reprove Thomassen's result on $C_5$-sub-free
graphs, and, for instance, gives that the chromatic threshold
is also 0 for Petersen-sub-free graphs.

The tool we use is a variation of VC-dimension. In the next section,
we show how the classical VC-theory is the key-tool for triangle-free
graphs, both giving an easy proof of the 1/3 threshold and
going slightly beyond.

\section{VC-dimension and Dense Triangle-Free Graphs}

Let $H=(V,E)$ be a hypergraph. A subset $X$ of $V$ is {\it shattered} by $H$
if for every subset $Y$ of $X$, there exists an edge $F$ of $H$ such that
$F\cap X=Y$. Introduced in \cite{Sau} and \cite{VC},
the {\it Vapnik-\v{C}ervonenkis dimension}
(or {\it VC-dimension}) of a hypergraph $H=(V,E)$ is the
maximum size of a subset shattered by $H$.

A {\it transversal} of $H$ is a set $T$ of vertices which intersects every edge
of $H$. The {\it transversality} of $H$, denoted by $\tau (H)$,
 is the minimum size of a transversal.
A  {\it fractional transversal} is a weight function
$w$ on the vertices of $H$ such that every edge has weight at least 1.
The {\it fractional transversality} $\tau ^*(H)$ of~$H$
is the minimum value of $w(V)$ taken over
all fractional transversal $w$ of $H$.
Since every transversal is a fractional transversal once interpreted as a
0,1 function, we clearly have that $\tau ^*(H)\leq \tau (H)$. However, there is no bound
on $\tau$ depending only on $\tau ^*$. Indeed, the hypergraph which
hyperedges are the $n$-element subsets of a $2n$-element  vertex set
have fractional transversality 2 and transversality $n+1$. The following
bound~\cite{HW} is one of the most useful application of VC-dimension.

\begin{lem} \label{transversal} Every hypergraph $H$ with VC-dimension $d$ satisfies
$$\tau (H)\leq 16d\tau ^*(H)\log (d\tau^*(H))$$
\end{lem}

The following result is a straightforward application of Lemma~\ref{transversal}. By {\it cube},
we mean the graph of the 3-dimensional cube.

\begin{theo} \label{tricube}
Every triangle-free graph $G$ with minimum degree $cn$ and without induced
cube satisfies $$\chi(G)\leq {192\log (6/c)\over c}$$
\end{theo}

\pro Consider a maximum cut $(X,Y)$ of $G$, i.e. a partition of the vertex
set which maximizes the number of edges between $X$ and $Y$.
Our goal is to prove that the chromatic number of the induced restriction
of $G$ on $X$ is at most  $96\log (6/c)/c$.
By maximality of the cut, every vertex in $X$ has as many neighbors in $Y$ as in
$X$. We now consider the hypergraph $H$ on vertex set $Y$
whose hyperedges are the neighborhoods in $Y$ of the vertices of $X$.
By the previous remark,
the size of every hyperedge of $H$ is at least $cn/2$. Hence the fractional transversality
of $H$ is at most $2/c$ (just observe that the constant weight function
valued $2/cn$ is a fractional transversal). The crucial point is the following: since $G$ is cube-free,
there is no vertices $x_1,x_2,x_3,x_4$ in $X$ and $y_1,y_2,y_3,y_4$ in $Y$
for which $x_i$ is not joined to $y_i$ for every $i=1,\dots ,4$ and $x_i$ is joined to
$y_j$ for all $i\neq j$. Interpreted in the hypergraph $H$, this means
that given a set $S$ of four vertices of $H$, one cannot find four hyperedges containing exactly
the four 3-elements subsets of $S$. Thus the VC-dimension of $H$ is at most three,
hence by Lemma~\ref{transversal} it has a transversal $T$ of size at most $96\log (6/c)/c$.
By definition of $H$, every vertex in $X$ has a neighbor in $T$, in other words,
the set of neighborhoods of the vertices of $T$ covers $X$. Since $G$ is triangle-free,
the neighborhood of a vertex is a stable set. In particular, the chromatic
number of the induced restriction of $G$ on $X$ is at most $96\log (6/c)/c$.
The same bound also holds for $Y$.~$\square$

~

Brandt proved in~\cite{BRA} that every maximal triangle-free graph with minimum degree
greater than $n/3$ has no induced cube.
We follow his argument in the proof of the next result.

\begin{theo}
Every triangle-free graph with minimum degree at least $n/3$ has chromatic number at most $1665$.
\end{theo}

\pro We first make $G$ maximal triangle-free by iteratively joining pairs
of vertices $x,y$ by an edge when their distance is at least 3. If $G$ does not contain
an induced cube, then from  Theorem~\ref{tricube} the chromatic number of
$G$ is bounded from above by $576\log (18)<1665$.

Now $G$ contains an induced cube consisting of the vertices $x_1,x_2,x_3,x_4$ and
$y_1,y_2,y_3,y_4$ as described in the proof of Theorem~\ref{tricube}. Since $x_i$ is not linked to $y_i$, their distance
is exactly two, hence there
exists a vertex $z_i$ joined to both $x_i$ and $y_i$. Observe that these vertices $z_i$
are pairwise distinct since $G$ is triangle-free. We denote by $C$ the graph
induced by these 12 vertices.

First observe that the maximum stable sets of $C$ have at most
five vertices, and that such a set is necessarily of the form
$x_i,y_i$ and all $z_j$ with $j\neq i$.

Now since $C$ has twelve vertices, the total number of edges incident
to $C$ is at least $4n$. If every vertex of $G$ has a neighbor
in $C$, the graph $G$ is covered by the union of twelve neighborhoods,
hence has chromatic number at most 12. Thus some vertex of $G$ has
no neighbor in $C$, which implies that some vertex has at least five
neighbors in $C$. We assume without loss of generality that this vertex is joined
to $x_4,y_4,z_1,z_2,z_3$. We call this new vertex $z_4$ and discard the previous
$z_4$. Now let $C'$ be the graph induced on $x_i,y_i,z_i$ for $i=1,\dots ,3$.
There are at least $3n$ edges incident to $C'$, and provided that the chromatic number of
$G$ is more than nine,
there exists a vertex which is not joined to $C'$, hence some vertex must be
joined to four vertices of $C'$. Again the unique possible choice is $x_i,y_i$ and all
$z_j$ with $j\neq i$. Let us say that there exists a vertex $z_1$ joined to $x_1,y_1,z_2,z_3$.
To avoid triangles, $z_1$ is not joined to $z_4$, so the vertices $z_1,z_2,z_3,z_4$
induce a cycle of length four. We still call $C$ this graph with the new vertices $z_1$ and $z_4$.
Now the maximum size of a stable set in the graph $C$ is four. Hence every vertex of $G$
sees exactly four neighbors in $C$, and thus $G$ has chromatic number at most twelve.~$\square$

~

A more careful analysis of the previous proof gives the following result, which breaks
the $n/3$ threshold.

\begin{theo}
Every triangle-free graph $G$ with $n$ vertices
 and minimum degree at least $n/3-b$, has chromatic number
at most  $\max\{10^5, 12(b+1)\}$.
\end{theo}

Let us remark that since triangle-free graphs may have
arbitrarily large chromatic number, so
for every constant $a$ there exist a constant $b$ and a
triangle-free graph $G$ with $n$ vertices and minimum degree
larger than $n/3-b$ such that $\chi(G)\ge a$.


\section{Paired VC-dimension}

Let $H=(V,E)$ be a hypergraph. A subset $F=\{F_1,F_2,\dots ,F_k\}$ of edges of $E$
forms a {\it complete Venn diagram} if for every subset $I\subseteq \{1,\dots ,k\}$,
there exists a vertex $v$ such that $v\in F_i$ for all $i\in I$, and $v\notin F_j$
for all $j\in  \{1,\dots ,k\}\setminus I$.

The {\it dual VC-dimension} of a hypergraph $H=(V,E)$ is the maximum size
of a subset of hyperedges which forms a complete Venn diagram. This is precisely
the VC-dimension of the dual hypergraph of $H$, obtained by exchanging
the roles of the vertices and the edges in the incidence bipartite graph of $H$.

Observe that the dual VC-dimension is invariant if we add to the edge set
of $H$ the complement of the edges of $H$. Indeed what is crucial
for this dimension is the bipartition $(e,\overline{e})$ generated by an edge $e$ rather than the
edge $e$ itself. This suggests the following extension of the definition of dual VC-dimension to arbitrary
pairs of subsets.

A {\it paired hypergraph} $P$ is a couple $(H,G)$ where $H=(V,E)$ is a hypergraph
and $G$ is a graph on vertex set $E$. Hence the edges of $G$ are pairs of hyperedges
of $H$. The {\it paired VC-dimension} of $H$ is the maximum
$d$ for which there exists $d$ edges $(A_i,B_i)_{i\in \{1,\dots ,d\}}$ of $G$ such that
for every subset $I$ of $\{1,\dots ,d\}$, there exists a vertex which belongs to
all $A_i$ where $i\in I$ and all $B_j$ where $j\in  \{1,\dots ,d\}\setminus I$.
By the {\it restriction of $P=(H,G)$ to $W\subseteq V$}, denoted by $P|_W$,
we mean paired hypergraph $P|_W=(H_W,G_W)$, where $H_W=(W,E_W)$,
$E_W=\{e\cup W:e\in E\}$, and $\{e',f'\}$ is an edge of $G_W$ if for some
edge $\{e,f\}$ of $G$ we have $e'=e\cap W$, $f'=f\cap W$.

Observe that if $H$ is a hypergraph, the paired VC-dimension of the set
of pairs $(e,\overline{e})$ is precisely its dual VC-dimension.

Our main result is an analogue of Lemma~\ref{transversal} for paired VC-dimension.

\begin{theo} \label{main}
There exists a function $f:[0,1]\times \N\ra \N$ such that every paired hypergraph $P=(H,G)$
with paired VC-dimension $d$ satisfies $\chi (G)\leq f(\tau ^*(H),d)$.
\end{theo}

\pro We let $H=(V,E)$. Observe first that duplicating any vertex $v$ into
$v_1,v_2$ (and replacing any occurrence of $v$ by $v_1,v_2$ in the edges of $E$)
does not affect $\chi$, $\tau ^*$ and $d$. Let us fix $c:=1/\tau^*$. One can duplicate
the vertices of $H$ accordingly to $\tau ^*$ in such a way that every edge
of $E$ has at least $cn$ vertices where $n:=|V|$.

Let us now fix some very small $\varepsilon >0$.

The first step of the proof is based on some density increase argument.
Given a subset $S$ of $V$, an edge $F$
of $E$ is {\it boosted} by $S$ if $${|F\cap S|\over |S|}\geq (1+\varepsilon){|F|\over |V|},$$ i.e.
its density in $S$ is larger than in $V$ by a factor at least $1+\varepsilon$.

A {\it $p$-booster} is a set
of nonempty subsets $S_1,\dots ,S_p$ of $V$ such that $E$ can be partitioned
into $E_0,\dots ,E_p$, where every edge of $E_i$, for $1\leq i\leq p$,
is boosted by $S_i$, and $E_0$ is a stable set of $G$. If some hyperedge
of $E$ is boosted by several $S_i$, we arbitrarily put it in some corresponding $E_i$.

Hence if such a $p$-booster exists, construct
the family of  $p$ paired hypergraphs
$P_1,\dots, P_p$, where $P_i=P|_{S_i}$,
and the stable set $E_{0}$. 
This process is iterated whenever one of the paired hypergraph
$P_1,\dots, P_p$ has in turn a $p$-booster. This gives a $(p+1)$-ary tree, called
{\it booster tree} where every internal
node corresponds to a $p$-booster, and which depth is bounded since the density
of every hyperedge is multiplied by $1+\varepsilon$ at each step and every hyperedge
starts with initial density at least $c$. Hence the size of this tree is bounded in terms
of $\varepsilon$, $p$ and $c$. Note that every leaf of this
tree is either a stable set or some paired hypergraph $P_i$ which has no $p$-booster.

We will slightly modify our booster tree in order to only manipulate hypergraphs
with hyperedges of the same size. This will simplify a bit our argument. Indeed, we can
assume that every hyperedge of $H=(V,E)$ has size exactly $cn$ since we can
arbitrarily delete some vertices from the hyperedges of $H$ with size more than $cn$.
Restricting hyperedges to fewer vertices can only increase the chromatic number
(since some hyperedges can be identified) and can only decrease
the paired VC-dimension. Similarly, during the construction of the booster tree,
we can reduce the density of the hyperedges of the nodes of depth $s$ to $(1+\varepsilon)^sc$.

The aforementioned construction of the booster tree can be done for
any value of $p$. For our proof, we tune $p$ to be equal to $2^{d}+1$.
To achieve our conclusion, we just have
to show that every leaf of the booster tree
has bounded chromatic number. This is obvious for the leaves which are
stable sets. We have to prove it for the leaves corresponding to some paired hypergraph $P_i$
which has no $(2^{d}+1)$-booster (and hence no $p$-booster for $p\leq 2^{d}+1$).

We now rename $P_i$ into $P=(H,G)$. We still assume that $H=(V,E)$ has $n$ vertices and
every hyperedge of $H$ has linear density $(1+\varepsilon)^sc$.
We still call $c$ this density $(1+\varepsilon)^sc$. Hence the hyperedges
have size $cn$. The paired
VC-dimension of $P$ is now at most $d$ (since we are dealing with a restriction of our original
hypergraph). Finally, we assume that $P$ has no $(2^{d}+1)$-booster.

The second step of our proof is to suppress the edges of $G$ joining two hyperedges of $H$ with
large intersection. The  {\it overlap} of an edge $(A_i,B_i)$ of $G$ is
the intersection $A_i\cap B_i$.
Let us set a positive constant  $\varepsilon ' < \varepsilon (c/9)^{2p}$.
We denote by $G'$ the subgraph of $G$ which consists
of edges of $G$ with overlap size more than $\varepsilon' n$. If $G'$ has a matching
of size more than $d/\varepsilon'$, there exists a vertex $v$ of $V$
which belongs to the overlap of $d+1$ edges of this matching. Hence, the paired VC-dimension
is at least $d+1$, a contradiction. Therefore the maximum size of a
matching is at most $d/\varepsilon'$, and thus the chromatic number of $G'$ is at most
$2d/\varepsilon'+1$ (since removing the vertices of a maximum matching leaves a stable set).
Delete all the edges of $G'$ from the graph $G$ so
that the remaining graph has only edges with small overlap and,
to simplify the notation, rename  the resulting graph as $G$.
In order to complete the proof it is
enough to show that the chromatic number of $G$ is bounded.

Observe first that if $(A_i,B_i)$ is an edge of $G$, since the overlap is at most
$\varepsilon' n$ and $A_i$ and $B_i$ have equal size, the density
of both $A_i$ and $B_i$ is at most $(1+\varepsilon')/2$. In particular $c<1/(1+\varepsilon)$
since $\varepsilon $ is small.

Since $P$ has no 0-booster, $G$ does not induce a stable set, we can find a pair $(A_1,B_1)$ in $G$.
By the definition of $G$ and $H$, both $A_1$ and $B_1$ have size at least $cn$. The overlap
of $A_1$ and $B_1$ is at most $\varepsilon' n$ which, let us recall,
is much smaller than  $\varepsilon (c/9)^{2p}<cn$.
Hence the subsets $A_1,B_1,C_1:=V\setminus (A_1\cup B_1)$
is a partition of $V$
up to $\varepsilon' n$ vertices so we can treat it as a real partition
in the remaining part of the argument.  By this we mean that in the density
estimates which follow we shall treat $A_1$ and $B_1$ as disjoint, keeping in
mind that in this way we can be off at most by a term of  $2\varepsilon' /c$
(which becomes $(2\varepsilon' /c)^i$ in the $i$th step of the procedure).

We now show the crucial argument - it will be repeated $d$ times to reach a contradiction.  The
subsets $A_1,B_1,C_1$ do not form a 3-booster, hence there exists an edge $(A_2,B_2)$ of $G$ such
that neither $A_2$ nor $B_2$ is boosted by $A_1$, $B_1$ or $C_1$. Observe that since $c<1/(1+\varepsilon)$,
$A_1$ is boosted by $A_1$. Hence the hyperedge $A_2$ is distinct from $A_1$. By the same argument, the
edges $(A_1,B_1)$ and $(A_2,B_2)$ of $G$ have no common endpoint.

By definition, the density of $A_2$ inside $A_1$, $B_1$ and $C_1$ is no more than $(1+\varepsilon)c$.
The key-point is that the density $\delta $ of $A_2$ in $A_1$ is close to $c$. Indeed, we have
$|A_1\cap A_2|=c\delta n$, and therefore the density of $A_1$ in $B_1\cup C_1$ is at least
$(c-c\delta)/(1-c)$. Since the density of $A_1$ in $B_1$ and $C_1$ is at most
$(1+\varepsilon)c$, the density of $A_1$ in $B_1\cup C_1$ is also at most
$(1+\varepsilon)c$. Thus $(1-\delta)/(1-c)\leq 1+\varepsilon$, and then $\delta> c+\varepsilon (c-1)$. Finally,
the density of $A_2$ in $A_1$ is roughly equal to $c$.

We now form four near disjoint subsets $A_1\cap A_2, A_1\cap B_2,A_2\cap B_1,A_2\cap B_2$
of size roughly $c^2n$ which, together with the complement of their union in $V$, are five subsets of $V$ not
forming a 5-booster. Hence there exists an edge $(A_3,B_3)$ of $G$ such that neither $A_3$ nor $B_3$
is boosted by these subsets. Again the edge $(A_3,B_3)$ has no common endpoint with $(A_1,B_1)$ or
$(A_2,B_2)$.  Moreover, the density of $A_3$ inside, say, $A_1\cap A_2$ is again roughly $c$, and
we can form eight subsets of $V$ formed by the intersections of the endpoints of $(A_1,B_1),(A_2,B_2)$ and
$(A_3,B_3)$. These eight subsets together with the leftover vertices do not form a 9-booster,
and an edge $(A_4,B_4)$ can be found.

This process extract a matching $(A_1,B_1),(A_2,B_2),\dots , (A_k,B_k)$ of $G$, with
$k$ at most $d+1$, since we use the fact that no $(2^d+1)$-booster exist. If $k=d+1$,
since the paired VC-dimension is equal to $d$, one of the intersections, say
$A_1\cap A_2\cap \dots \cap A_k$, is empty. But since this intersection has
size roughly $c^{d+1}n$, this mean that $n=O(c^{-d-1})$ and hence the size
of $G$ (and in particular its chromatic number) is bounded in terms of $c$ and $d$.
If $k<d+1$, the process has stopped on $(A_1,B_1),(A_2,B_2),\dots , (A_k,B_k)$
while no $2^k+1$-booster exists. The sole reason for this is that one of the intersection
is empty, again implying that $n=O(c^{-k})$ again a bound in terms of $c$ and $d$.

This completes the proof of Theorem~\ref{main}.~$\square$

~

Let us make two remarks concerning the above result. From
the proof it follows that if the paired hypergraph $P$ has
fixed $\tau ^*$ and $G$ has arbitrarily high chromatic number, then the paired VC-dimension is unbounded, which is witnessed by pairs $(A_i,B_i)$
forming a matching of $G$.

Let us now introduce two natural paired-hypergraphs associated to a graph $G=(V,E)$.

\begin{itemize}
\item The {\it neighborhood paired hypergraph} $P_N=(H_N,G_N)$ in  which the vertices
of $H_N$ are the vertices of $G$, the hyperedges of $H_N$ are the neighborhoods of
the vertices of $G$, and where two neighborhoods $N(v)$ and $N(w)$ are neighbors in
$G_N$ if $vw$ is an edge of $G$.

\item The {\it stable paired hypergraph} $P_S=(H_S,G_S)$ in which the vertices
of $H_S$ are the stable sets of $G$, the hyperedges $S_v$ of $H_S$ are the stable sets
containing a given vertex $v$ of $G$, and where two hyperedges $S_v$ and $S_w$ are neighbors in
$G_S$ if $vw$ is an edge of $G$.

\end{itemize}

As Pierre Charbit observed, the fractional transversality of the stable paired hypergraph is a classical graph parameter.
To see this, note that the bipartite incidence graph of $H_S$ is isomorphic to
the incidence bipartite graph of vertices versus stable sets. Hence $\tau ^*$ is the
minimum sum of a weight function on stable sets such that for every vertex $v$, the
sum of the weights of stable sets containing $v$ is at least 1. This is exactly
the {\it fractional chromatic number} of $G$, denoted by $\chi ^*(G)$. Hence
Theorem~\ref{main} can be reformulated as:

\begin{theo} \label{stable}
Every graph $G$ satisfies $\chi (G)\leq f(\chi ^*(G),d)$ where $d$ is the paired VC-dimension
of its stable paired hypergraph.
\end{theo}

In this form, this result is not really interesting since $d$
can be expressed in terms of the stability of $G$, and certainly $\chi (G)$
is bounded in terms of $\chi ^*(G)$ and $\alpha (G)$.
It would be interesting to modify the stable paired hypergraph
construction (maybe by selecting some special stable sets)
to possibly derive some non-trivial inequality relating $\chi (G)$ and $\chi^* (G)$.

\section{Chromatic Threshold for Homomorphisms.}

A graph $G=(V,E)$ is a {\it near bipartite graph} if $G$ is triangle-free and
admits a vertex partition $V=V_1\cup V_2$ such that $V_1$ induces a stable set
and $V_2$ induces a graph with maximum degree one (i.e. a partial matching).

The simplest examples of near bipartite graphs which are not bipartite are odd
cycles. However some non-trivial graphs like the Petersen graph are
also near bipartite since removing any stable set of size four leaves an
induced  matching of three edges.

A construction of a family of universal near bipartite graph (i.e. containing every
near bipartite graph) can be obtained as follows. Given a graph $G$ with connected
components $C_1,\dots ,C_m$, we form a new
graph by adding, for every choice of vertices $v_1,\dots ,v_m$ in these components,
a new vertex linked exactly to these vertices. This construction was introduced
by Zykov~\cite{Zyk} in order to construct triangle-free graphs with arbitrarily
large chromatic number. We call {\it Zykov $k$-matching} the result of this construction
when the original graph $G$ is a matching on $k$ edges. We denote it by $Z_k$.
It is routine to check that the class of graphs $Z_k$ contain
all near bipartite graphs.

The following result is a straightforward consequence of Theorem~\ref{main}.

\begin{theo} \label{homfree}
The class of $H$-hom-free graphs has chromatic threshold 0 if and only if $H$
admits a homomorphism into a near bipartite graph.
\end{theo}

\pro
We prove first that if $H$ is not homomorphic to a near bipartite graph, its chromatic threshold
is at least 1/3. For this, consider a Borsuk-Hajnal graph $BH$ with near odd girth larger than the number
of vertices of $H$. Assume now for contradiction that there exists a homomorphism from $H$
to $BH$. Let $G$ be the graph induced by the image of $H$. Denote by $S$ the stable set of $G$
consisting of vertices of $G$ which are caps of $BH$. Since the near odd girth is large enough,
$G\setminus S$ is a bipartite graph. We denote the partite sets of the connected components
by $(A_1,B_1), (A_2,B_2), \dots ,(A_k,B_k)$. Since the near odd girth is large enough, no vertex of $S$
can be a neighbor of both a vertex of $A_i$ and a vertex of $B_i$. Hence contracting
$(A_i,B_i)$ to an edge $\{a_i,b_i\}$ for every $i$ does not create a triangle,
and therefore is a homomorphism from $G$ (hence from $H$) into a near bipartite graph, a contradiction.

Assume now that $H$ has a homomorphism into some near bipartite graph, hence
into some $Z_k$.
Our goal is to prove that the class of $Z_k$-hom-free graphs has chromatic threshold 0. Consider for
this a $Z_k$-hom-free graph $G=(V,E)$ with $n$ vertices and minimum degree $\delta \geq c n$
for some fixed $c>0$. Our goal is to show that $\chi (G)$ is bounded in terms of $c$.
To see this, we consider the neighborhood paired hypergraph $P_N=(H_N,G_N)$
of $G$.

The first remark is that $\tau^*(H)$ is at most $1/c$ since $\delta$
is at least $cn$. Thus we just have to
bound the paired VC-dimension of $P_N$ by some fixed constant to conclude
the statement from Theorem~\ref{main}.
Observe that $G$ is triangle-free since otherwise
$Z_k$, which is 3-chromatic, would be homomorphically mapped into the triangle.
Hence $P_N$ is a paired hypergraph in which the paired subsets do not
overlap. Assume for contradiction that
its paired VC-dimension is at least $k$. There exists a set of $k$ edges $(a_1,b_1), \dots ,(a_k,b_k)$
of $G$ corresponding to $k$ pairs of neighborhoods $(A_1,B_1),\dots , (A_k,B_k)$ with
paired VC-dimension $k$.
We claim that $Z_k$ maps into $G$ simply by identifying the $k$-matching of $Z_k$ with the edges
$(a_1,b_1), \dots ,(a_k,b_k)$. Note also that
 one can easily map to $G$ the other vertices of $Z_k$ which are joined to exactly
one vertex per edge of the matching. For instance, a vertex of
$G$ joined to $a_1,b_2,b_3,\dots ,B_k$, should be mapped
into a vertex from the intersection $A_1\cap B_2\cap B_3 \dots \cap B_k$,
 which, by the definition of paired VC-dimension, is nonempty. This contradicts
the fact that the paired VC-dimension is at least $k$.~$\square$

\section{Chromatic Threshold for Subgraphs.}

In the previous section, we proved that the class of graphs which does not
contain a homomorphic pentagon or Petersen graph has chromatic
threshold equal to zero. In the case of subgraphs, Thomassen~\cite{CT2}
proved the following result.

\begin{theo} \label{pentagon}
The class of $C_{2k+1}$-sub-free graphs has chromatic threshold 0, provided $k\ge2$.
\end{theo}

Again, paired VC-dimension is a tailor-made tool for this kind of questions. For instance,
we can extend Thomassen's result to all near bipartite graphs
rather than just odd cycles.
Unfortunately, unlike the homomorphism case, we could not characterize the graphs $H$ for
which the chromatic threshold is equal to zero. We discuss a possible characterization
in the final section.

\begin{theo} \label{subfree}
If $H$ is a near bipartite graph, the class of $H$-sub-free graphs has chromatic threshold 0.
\end{theo}

\pro
It suffices to prove the result when $H$ is isomorphic to some $Z_k$. Consider for
this a $Z_k$-sub-free graph $G=(V,E)$ with $n$ vertices and minimum degree $\delta \geq c n$,
for some fixed $c>0$. Assume that $G$ is chosen with arbitrarily high chromatic number, and
consider a maximum cut $V_1,V_2$ of $G$. Without loss of generality,
assume that $G_1$, the induced restriction of $G$
to $V_1$, has arbitrarily high chromatic number. Consider the set $V_2\choose 2^k$ of all
subsets of $V_2$ of size $2^k$, and form the following paired hypergraph $P=(H,G')$:

\begin{itemize}
\item the vertex set of $H$ is $V_2\choose 2^k$;
\item each hyperedge $e_{v_1}$ of $H$ corresponds to some vertex $v_1$ of $V_1$ and
contains all the $2^k$-element subsets of $V_2$ contained in the neighborhood
of $v_1$;
\item every edge $v_1v_2$ of $V_1$ corresponds to the edge $e_{v_1}e_{v_2}$
of $G'$.
\end{itemize}

In other words, we simply lift the graph on $V_1$ to the set of neighborhoods in $V_2\choose 2^k$.
Observe that the size of a hyperedge is still linear, although its density decreased to $c^{2^k}$.
Hence, since the chromatic number of $G$ is arbitrarily large, one can find a matching
$(A_1,B_1),(A_2,B_2),\dots ,(A_k,B_k)$ of $G'$ such that all the $2^k$ possible intersections are
nonempty. It can be the case that the same vertex of $H$ is used to provide the nonempty
intersection (as in the proof of  Theorem~\ref{main}).
But since a vertex of $H$ corresponds to a subset of size $2^k$ of $V_2$, one can provide
nonempty intersection using $2^k$ distinct vertices.

To conclude, note that this matching $(A_1,B_1),(A_2,B_2),\dots ,(A_k,B_k)$ corresponds to
a matching $(a_1,b_1),(a_2,b_2),\dots ,(a_k,b_k)$ of $V_1$, which can be extended to
a Zykov $k$-matching to reach a contradiction.~$\square$

\section{Conclusion and Open Problems.}

An obvious open question is the behavior of the function
$f$ in Theorem~\ref{main}. We did not make any attempt to
optimize it in our proof. We do not even know
if polynomial bound is achievable.

Concerning extremal graphs, one of the main question left
open in this paper is the characterization of graphs $H$
such that the class of $H$-sub-free graphs has chromatic threshold
0. A {\it near acyclic}
graph $G$ is a graph admitting a stable set $S$ such that
$G-S$ is a forest and such that every odd cycle has at least
two vertices in $S$.

Following a construction of Hell and
Ne{\v s}et{\v r}il, given a graph $G$, there exists a graph $H_G$
which is homomorphic to $G$, has same chromatic number as $G$,
and has arbitrary high girth. We use this to provide a graph
with high
chromatic number,
which all small induced restrictions are near acyclic,
and has minimum degree close to 1/3. Start
with some Borsuk graph $B$. Then, apply the Hell and
Ne{\v s}et{\v r}il method to form the graph $H_B$. Finally,
use this graph $H_B$ instead of $B$ to form a Borsuk-Hajnal
graph. We conjecture that this construction is generic.

\begin{conj}
The class of $H$-sub-free graphs has chromatic threshold 0 if and only if $H$
is a near acyclic graph.
\end{conj}

Again this would give a gap between 0 and 1/3.

Our last question is an attempt to generalize triangle-free
graphs. A graph is {\it locally bipartite} if the neighborhood
of every vertex induces a bipartite graph. We conjecture the following.

\begin{conj}
The class of locally bipartite graphs has chromatic threshold 1/2.
\end{conj}

Let us sketch the construction showing that degree close to 1/2
is achievable by some locally bipartite graph. Again the construction
is topological. Consider the $d$-dimensional sphere. Given a point $x$
of the sphere, the {\it $(x,\delta)$-bicap} is the set of points of the sphere
which distance to $x$ is at most $\pi/2-\delta$ or at least $\pi/2+\delta$.
Now form a graph with parameters $(d,\varepsilon,\delta)$  as follows:

\begin{itemize}
\item start with $B$, a Borsuk graph  $Bor(d,\varepsilon)$;
\item add a stable set $S$ which vertices correspond to all $(x,\delta)$-bicaps;
\item join $x\in B$ to $c\in S$ if the point $x$ belongs to the bicap $c$.
\end{itemize}

This graph has high chromatic number due to the Borsuk graph. Since
a bicap nearly cover the sphere (provided that $\delta$ is small
enough), the bipartite graph $B,S$ is nearly complete, hence
giving minimum degree close to 1/2. Finally, the neighborhood
of vertices are bipartite, provided that $\delta <\varepsilon$.

~

We gratefully thank Pierre Charbit and Jarik Ne{\v s}et{\v r}il for helpful comments.

This work was done while the authors were visiting IPAM.

\end{document}